\documentclass{article}
\usepackage{graphics}
\usepackage{tikz,tikz-qtree,tikz-qtree-compat,adjustbox}
\usepackage{pst-plot}
\usepackage{etoolbox}

\title{Rarity of the infinite chains in the tree of numerical semigroups}

\author{Maria Bras-Amorós, Mariana Rosas Ribeiro}

\usepackage{amsmath,amsthm}
\newtheorem{theorem}{Theorem}
\newtheorem{lemma}[theorem]{Lemma}

\begin{document}

\maketitle

\begin{abstract}
We prove that, for each fixed genus, the portion of semigroups of that genus belonging to infinite chains in the semigroup tree approaches $0$ as the genus grows to infinite. This means that most numerical semigroups have a finite number of descendants in the semigroup tree. This problem has been open since 2009.
  \end{abstract}

\section{Introduction}

A numerical semigroup is a cofinite submonoid of the nonnegative integers.
The nonnegative integers not belonging to a semigroup are called its gaps, and the number of gaps is its genus.
In the last years there have been many efforts to count the number of numerical semigroups of each given genus, either theoretically \cite{Br:Bounds,Elizalde,Zhao,BlancoGarciaPuerto,Kaplan,Zhai,ODorney,EliahouRamirezAlfonsin,BernardiniTorres,EliahouFromentin:gapsets,EliahouFromentin:gapsetsm,Zhu} and computationally \cite{Br:Fibonacci,FromentinHivert,seeds1,rgd,seeds2,DelgadoEliahouFromentin}. One main tool is that of the tree of numerical semigroups. The root of this tree is the trivial numerical semigroup of genus $0$, and at each depth of the tree there are the semigroups whose genus equals the depth. The parent of each non-trivial semigroup in the tree is the semigroup obtained by adding to the semigroup its largest gap, the so-called Frobenius number of the semigroup. This way, if one can explore the whole tree of semigroups, computing the number of semigroups of a given genus is equivalent to counting the number of nodes of the tree at the given depth.

\begin{figure}
%
%
%
%

\providecommand\circledcolorednumb{}\renewcommand\circledcolorednumb[2]{\resizebox{0.065404\textwidth}{!}{\tikz[baseline=(char.center)]{\node[shape = circle,draw, inner sep = 2pt,fill=#1](char)    {\phantom{00}};\node[anchor=center] at (char.center) {\makebox(0,0){\large{#2}}};}}}
\robustify{\circledcolorednumb}
\providecommand\nongap{}\renewcommand\nongap[1]{\circledcolorednumb{gray!40}{#1}}
\providecommand\gap{}\renewcommand\gap[1]{\circledcolorednumb{black!05}{\phantom{#1}}}
\providecommand\generator{}\renewcommand\generator[1]{\circledcolorednumb{gray}{#1}}
\providecommand\seed{}\renewcommand\seed[1]{\circledcolorednumb{black!30}{#1}}
\providecommand\nonseed{}\renewcommand\nonseed[1]{\circledcolorednumb{black!05}{#1}}
\providecommand\dotscircles{}\renewcommand\dotscircles{\resizebox{0.065404\textwidth}{!}{\dots}}
\providecommand\gapingapset{}\renewcommand\gapingapset[1]{\circledcolorednumb{gray!50}{#1}}
\providecommand\nongapingapset{}\renewcommand\nongapingapset[1]{\phantom{\gapingapset{#1}}}
\adjustbox{max width=\textwidth,max height=.9\textheight}{\begin{tikzpicture}[grow'=right, sibling distance=6.000000mm]\tikzset{level 1/.style={level distance=21.000000cm}}\tikzset{level 2/.style={level distance=26.250000cm}}\tikzset{level 3/.style={level distance=34.125000cm}}\tikzset{level 4/.style={level distance=42.000000cm}}\tikzset{level 5/.style={level distance=52.500000cm}}\tikzset{level 6/.style={level distance=53.550000cm}}\tikzset{level 7+/.style={level distance=56.700000cm}}\Tree[.{}  \edge [black,thick]; [.{}  \edge [black,thick]; [.{}  \edge [black,thick]; [.{}  \edge [black,thick]; [.{}  \edge [black,thick]; [.{}  \edge [black,thick]; [.{}  \edge [black,thick]; [.{}  \edge [black,thick]; [.{}  \edge [black,thick]; [.{}  \edge [black,thick]; [.{}  \edge [black,thick]; [.{} ] \edge [black,thick]; [.{} ] \edge [gray!50]; [.{} ] \edge [gray!50]; [.{} ] \edge [gray!50]; [.{} ] \edge [gray!50]; [.{} ] \edge [gray!50]; [.{} ] \edge [gray!50]; [.{} ] \edge [gray!50]; [.{} ] \edge [gray!50]; [.{} ] \edge [gray!50]; [.{} ]] \edge [black,thick]; [.{}  \edge [black,thick]; [.{} ] \edge [black,thick]; [.{} ] \edge [gray!50]; [.{} ] \edge [gray!50]; [.{} ] \edge [gray!50]; [.{} ] \edge [gray!50]; [.{} ] \edge [gray!50]; [.{} ] \edge [gray!50]; [.{} ] \edge [gray!50]; [.{} ]] \edge [gray!50]; [.{}  \edge [gray!50]; [.{} ] \edge [gray!50]; [.{} ] \edge [gray!50]; [.{} ] \edge [gray!50]; [.{} ] \edge [gray!50]; [.{} ] \edge [gray!50]; [.{} ] \edge [gray!50]; [.{} ]] \edge [gray!50]; [.{}  \edge [gray!50]; [.{} ] \edge [gray!50]; [.{} ] \edge [gray!50]; [.{} ] \edge [gray!50]; [.{} ] \edge [gray!50]; [.{} ] \edge [gray!50]; [.{} ]] \edge [gray!50]; [.{}  \edge [gray!50]; [.{} ] \edge [gray!50]; [.{} ] \edge [gray!50]; [.{} ] \edge [gray!50]; [.{} ] \edge [gray!50]; [.{} ]] \edge [gray!50]; [.{}  \edge [gray!50]; [.{} ] \edge [gray!50]; [.{} ] \edge [gray!50]; [.{} ] \edge [gray!50]; [.{} ]] \edge [gray!50]; [.{}  \edge [gray!50]; [.{} ] \edge [gray!50]; [.{} ] \edge [gray!50]; [.{} ]] \edge [gray!50]; [.{}  \edge [gray!50]; [.{} ] \edge [gray!50]; [.{} ]] \edge [gray!50]; [.{}  \edge [gray!50]; [.{} ]] \edge [gray!50]; [.{} ]] \edge [black,thick]; [.{}  \edge [black,thick]; [.{}  \edge [black,thick]; [.{} ] \edge [black,thick]; [.{} ] \edge [gray!50]; [.{} ] \edge [gray!50]; [.{} ] \edge [gray!50]; [.{} ] \edge [gray!50]; [.{} ] \edge [gray!50]; [.{} ] \edge [gray!50]; [.{} ]] \edge [gray!50]; [.{}  \edge [gray!50]; [.{} ] \edge [gray!50]; [.{} ] \edge [gray!50]; [.{} ] \edge [gray!50]; [.{} ] \edge [gray!50]; [.{} ] \edge [gray!50]; [.{} ] \edge [gray!50]; [.{} ]] \edge [gray!50]; [.{}  \edge [gray!50]; [.{} ] \edge [gray!50]; [.{} ] \edge [gray!50]; [.{} ] \edge [gray!50]; [.{} ] \edge [gray!50]; [.{} ]] \edge [gray!50]; [.{}  \edge [gray!50]; [.{} ] \edge [gray!50]; [.{} ] \edge [gray!50]; [.{} ] \edge [gray!50]; [.{} ]] \edge [gray!50]; [.{}  \edge [gray!50]; [.{} ] \edge [gray!50]; [.{} ] \edge [gray!50]; [.{} ]] \edge [gray!50]; [.{}  \edge [gray!50]; [.{} ] \edge [gray!50]; [.{} ]] \edge [gray!50]; [.{}  \edge [gray!50]; [.{} ]] \edge [gray!50]; [.{} ]] \edge [gray!50]; [.{}  \edge [gray!50]; [.{}  \edge [gray!50]; [.{} ] \edge [gray!50]; [.{} ] \edge [gray!50]; [.{} ] \edge [gray!50]; [.{} ] \edge [gray!50]; [.{} ] \edge [gray!50]; [.{} ]] \edge [gray!50]; [.{}  \edge [gray!50]; [.{} ] \edge [gray!50]; [.{} ] \edge [gray!50]; [.{} ] \edge [gray!50]; [.{} ]] \edge [gray!50]; [.{}  \edge [gray!50]; [.{} ] \edge [gray!50]; [.{} ] \edge [gray!50]; [.{} ]] \edge [gray!50]; [.{}  \edge [gray!50]; [.{} ] \edge [gray!50]; [.{} ]] \edge [gray!50]; [.{}  \edge [gray!50]; [.{} ]] \edge [gray!50]; [.{} ]] \edge [gray!50]; [.{}  \edge [gray!50]; [.{}  \edge [gray!50]; [.{} ] \edge [gray!50]; [.{} ] \edge [gray!50]; [.{} ] \edge [gray!50]; [.{} ]] \edge [gray!50]; [.{}  \edge [gray!50]; [.{} ] \edge [gray!50]; [.{} ] \edge [gray!50]; [.{} ]] \edge [gray!50]; [.{}  \edge [gray!50]; [.{} ] \edge [gray!50]; [.{} ]] \edge [gray!50]; [.{}  \edge [gray!50]; [.{} ]] \edge [gray!50]; [.{} ]] \edge [gray!50]; [.{}  \edge [gray!50]; [.{}  \edge [gray!50]; [.{} ] \edge [gray!50]; [.{} ] \edge [gray!50]; [.{} ]] \edge [gray!50]; [.{}  \edge [gray!50]; [.{} ] \edge [gray!50]; [.{} ]] \edge [gray!50]; [.{}  \edge [gray!50]; [.{} ]] \edge [gray!50]; [.{} ]] \edge [gray!50]; [.{}  \edge [gray!50]; [.{}  \edge [gray!50]; [.{} ] \edge [gray!50]; [.{} ]] \edge [gray!50]; [.{}  \edge [gray!50]; [.{} ]] \edge [gray!50]; [.{} ]] \edge [gray!50]; [.{}  \edge [gray!50]; [.{}  \edge [gray!50]; [.{} ]] \edge [gray!50]; [.{} ]] \edge [gray!50]; [.{}  \edge [gray!50]; [.{} ]] \edge [gray!50]; [.{} ]] \edge [black,thick]; [.{}  \edge [black,thick]; [.{}  \edge [black,thick]; [.{}  \edge [black,thick]; [.{} ] \edge [black,thick]; [.{} ] \edge [gray!50]; [.{} ] \edge [gray!50]; [.{} ] \edge [gray!50]; [.{} ] \edge [gray!50]; [.{} ] \edge [gray!50]; [.{} ]] \edge [gray!50]; [.{}  \edge [gray!50]; [.{} ] \edge [gray!50]; [.{} ] \edge [gray!50]; [.{} ] \edge [gray!50]; [.{} ] \edge [gray!50]; [.{} ] \edge [gray!50]; [.{} ]] \edge [gray!50]; [.{}  \edge [gray!50]; [.{} ] \edge [gray!50]; [.{} ] \edge [gray!50]; [.{} ] \edge [gray!50]; [.{} ] \edge [gray!50]; [.{} ]] \edge [gray!50]; [.{}  \edge [gray!50]; [.{} ] \edge [gray!50]; [.{} ] \edge [gray!50]; [.{} ]] \edge [gray!50]; [.{}  \edge [gray!50]; [.{} ] \edge [gray!50]; [.{} ]] \edge [gray!50]; [.{}  \edge [gray!50]; [.{} ]] \edge [gray!50]; [.{} ]] \edge [black,thick]; [.{}  \edge [black,thick]; [.{}  \edge [black,thick]; [.{} ] \edge [gray!50]; [.{} ] \edge [gray!50]; [.{} ] \edge [gray!50]; [.{} ] \edge [gray!50]; [.{} ]] \edge [black,thick]; [.{}  \edge [black,thick]; [.{} ] \edge [black,thick]; [.{} ] \edge [gray!50]; [.{} ] \edge [gray!50]; [.{} ] \edge [gray!50]; [.{} ]] \edge [gray!50]; [.{}  \edge [gray!50]; [.{} ] \edge [gray!50]; [.{} ] \edge [gray!50]; [.{} ]] \edge [gray!50]; [.{}  \edge [gray!50]; [.{} ] \edge [gray!50]; [.{} ]] \edge [gray!50]; [.{}  \edge [gray!50]; [.{} ]] \edge [gray!50]; [.{} ]] \edge [gray!50]; [.{}  \edge [gray!50]; [.{}  \edge [gray!50]; [.{} ] \edge [gray!50]; [.{} ] \edge [gray!50]; [.{} ]] \edge [gray!50]; [.{}  \edge [gray!50]; [.{} ] \edge [gray!50]; [.{} ]] \edge [gray!50]; [.{}  \edge [gray!50]; [.{} ]] \edge [gray!50]; [.{} ]] \edge [gray!50]; [.{}  \edge [gray!50]; [.{}  \edge [gray!50]; [.{} ] \edge [gray!50]; [.{} ]] \edge [gray!50]; [.{}  \edge [gray!50]; [.{} ]] \edge [gray!50]; [.{} ]] \edge [gray!50]; [.{}  \edge [gray!50]; [.{}  \edge [gray!50]; [.{} ]] \edge [gray!50]; [.{} ]] \edge [gray!50]; [.{}  \edge [gray!50]; [.{} ]] \edge [gray!50]; [.{} ]] \edge [gray!50]; [.{}  \edge [gray!50]; [.{}  \edge [gray!50]; [.{}  \edge [gray!50]; [.{} ] \edge [gray!50]; [.{} ] \edge [gray!50]; [.{} ] \edge [gray!50]; [.{} ] \edge [gray!50]; [.{} ]] \edge [gray!50]; [.{}  \edge [gray!50]; [.{} ] \edge [gray!50]; [.{} ] \edge [gray!50]; [.{} ]] \edge [gray!50]; [.{}  \edge [gray!50]; [.{} ] \edge [gray!50]; [.{} ]] \edge [gray!50]; [.{}  \edge [gray!50]; [.{} ]] \edge [gray!50]; [.{} ]] \edge [gray!50]; [.{}  \edge [gray!50]; [.{}  \edge [gray!50]; [.{} ] \edge [gray!50]; [.{} ] \edge [gray!50]; [.{} ]] \edge [gray!50]; [.{}  \edge [gray!50]; [.{} ]] \edge [gray!50]; [.{} ]] \edge [gray!50]; [.{}  \edge [gray!50]; [.{}  \edge [gray!50]; [.{} ]] \edge [gray!50]; [.{} ]] \edge [gray!50]; [.{}  \edge [gray!50]; [.{} ]] \edge [gray!50]; [.{} ]] \edge [gray!50]; [.{}  \edge [gray!50]; [.{}  \edge [gray!50]; [.{}  \edge [gray!50]; [.{} ] \edge [gray!50]; [.{} ] \edge [gray!50]; [.{} ]] \edge [gray!50]; [.{}  \edge [gray!50]; [.{} ]] \edge [gray!50]; [.{} ]] \edge [gray!50]; [.{}  \edge [gray!50]; [.{}  \edge [gray!50]; [.{} ]] \edge [gray!50]; [.{} ]] \edge [gray!50]; [.{}  \edge [gray!50]; [.{} ]] \edge [gray!50]; [.{} ]] \edge [gray!50]; [.{}  \edge [gray!50]; [.{}  \edge [gray!50]; [.{}  \edge [gray!50]; [.{} ]] \edge [gray!50]; [.{} ]] \edge [gray!50]; [.{}  \edge [gray!50]; [.{} ]] \edge [gray!50]; [.{} ]] \edge [gray!50]; [.{}  \edge [gray!50]; [.{}  \edge [gray!50]; [.{} ]] \edge [gray!50]; [.{} ]] \edge [gray!50]; [.{}  \edge [gray!50]; [.{} ]] \edge [gray!50]; [.{} ]] \edge [black,thick]; [.{}  \edge [black,thick]; [.{}  \edge [black,thick]; [.{}  \edge [black,thick]; [.{}  \edge [black,thick]; [.{} ] \edge [gray!50]; [.{} ] \edge [gray!50]; [.{} ] \edge [gray!50]; [.{} ] \edge [gray!50]; [.{} ] \edge [gray!50]; [.{} ]] \edge [gray!50]; [.{}  \edge [gray!50]; [.{} ] \edge [gray!50]; [.{} ] \edge [gray!50]; [.{} ] \edge [gray!50]; [.{} ] \edge [gray!50]; [.{} ]] \edge [gray!50]; [.{}  \edge [gray!50]; [.{} ] \edge [gray!50]; [.{} ] \edge [gray!50]; [.{} ] \edge [gray!50]; [.{} ]] \edge [gray!50]; [.{}  \edge [gray!50]; [.{} ] \edge [gray!50]; [.{} ]] \edge [gray!50]; [.{}  \edge [gray!50]; [.{} ]] \edge [gray!50]; [.{} ]] \edge [gray!50]; [.{}  \edge [gray!50]; [.{}  \edge [gray!50]; [.{} ] \edge [gray!50]; [.{} ] \edge [gray!50]; [.{} ] \edge [gray!50]; [.{} ] \edge [gray!50]; [.{} ]] \edge [gray!50]; [.{}  \edge [gray!50]; [.{} ] \edge [gray!50]; [.{} ] \edge [gray!50]; [.{} ]] \edge [gray!50]; [.{}  \edge [gray!50]; [.{} ] \edge [gray!50]; [.{} ]] \edge [gray!50]; [.{}  \edge [gray!50]; [.{} ]] \edge [gray!50]; [.{} ]] \edge [gray!50]; [.{}  \edge [gray!50]; [.{}  \edge [gray!50]; [.{} ] \edge [gray!50]; [.{} ] \edge [gray!50]; [.{} ]] \edge [gray!50]; [.{}  \edge [gray!50]; [.{} ] \edge [gray!50]; [.{} ]] \edge [gray!50]; [.{}  \edge [gray!50]; [.{} ]] \edge [gray!50]; [.{} ]] \edge [gray!50]; [.{}  \edge [gray!50]; [.{}  \edge [gray!50]; [.{} ]] \edge [gray!50]; [.{} ]] \edge [gray!50]; [.{}  \edge [gray!50]; [.{} ]] \edge [gray!50]; [.{} ]] \edge [gray!50]; [.{}  \edge [gray!50]; [.{}  \edge [gray!50]; [.{}  \edge [gray!50]; [.{} ] \edge [gray!50]; [.{} ] \edge [gray!50]; [.{} ] \edge [gray!50]; [.{} ]] \edge [gray!50]; [.{}  \edge [gray!50]; [.{} ] \edge [gray!50]; [.{} ] \edge [gray!50]; [.{} ]] \edge [gray!50]; [.{}  \edge [gray!50]; [.{} ]] \edge [gray!50]; [.{} ]] \edge [gray!50]; [.{}  \edge [gray!50]; [.{}  \edge [gray!50]; [.{} ] \edge [gray!50]; [.{} ] \edge [gray!50]; [.{} ]] \edge [gray!50]; [.{}  \edge [gray!50]; [.{} ] \edge [gray!50]; [.{} ]] \edge [gray!50]; [.{}  \edge [gray!50]; [.{} ]] \edge [gray!50]; [.{} ]] \edge [gray!50]; [.{}  \edge [gray!50]; [.{}  \edge [gray!50]; [.{} ]] \edge [gray!50]; [.{} ]] \edge [gray!50]; [.{}  \edge [gray!50]; [.{} ]] \edge [gray!50]; [.{} ]] \edge [gray!50]; [.{}  \edge [gray!50]; [.{}  \edge [gray!50]; [.{}  \edge [gray!50]; [.{} ]] \edge [gray!50]; [.{} ]] \edge [gray!50]; [.{}  \edge [gray!50]; [.{} ]] \edge [gray!50]; [.{} ]] \edge [gray!50]; [.{}  \edge [gray!50]; [.{}  \edge [gray!50]; [.{} ]] \edge [gray!50]; [.{} ]] \edge [gray!50]; [.{}  \edge [gray!50]; [.{} ]] \edge [gray!50]; [.{} ]] \edge [gray!50]; [.{}  \edge [gray!50]; [.{}  \edge [gray!50]; [.{}  \edge [gray!50]; [.{}  \edge [gray!50]; [.{} ] \edge [gray!50]; [.{} ] \edge [gray!50]; [.{} ] \edge [gray!50]; [.{} ]] \edge [gray!50]; [.{}  \edge [gray!50]; [.{} ] \edge [gray!50]; [.{} ]] \edge [gray!50]; [.{}  \edge [gray!50]; [.{} ]] \edge [gray!50]; [.{} ]] \edge [gray!50]; [.{}  \edge [gray!50]; [.{}  \edge [gray!50]; [.{} ] \edge [gray!50]; [.{} ]] \edge [gray!50]; [.{} ]] \edge [gray!50]; [.{}  \edge [gray!50]; [.{} ]] \edge [gray!50]; [.{} ]] \edge [gray!50]; [.{}  \edge [gray!50]; [.{}  \edge [gray!50]; [.{}  \edge [gray!50]; [.{} ]] \edge [gray!50]; [.{} ]] \edge [gray!50]; [.{} ]] \edge [gray!50]; [.{}  \edge [gray!50]; [.{} ]] \edge [gray!50]; [.{} ]] \edge [gray!50]; [.{}  \edge [gray!50]; [.{}  \edge [gray!50]; [.{}  \edge [gray!50]; [.{}  \edge [gray!50]; [.{} ] \edge [gray!50]; [.{} ]] \edge [gray!50]; [.{} ]] \edge [gray!50]; [.{} ]] \edge [gray!50]; [.{}  \edge [gray!50]; [.{} ]] \edge [gray!50]; [.{} ]] \edge [gray!50]; [.{}  \edge [gray!50]; [.{}  \edge [gray!50]; [.{} ]] \edge [gray!50]; [.{} ]] \edge [gray!50]; [.{}  \edge [gray!50]; [.{} ]] \edge [gray!50]; [.{} ]] \edge [black,thick]; [.{}  \edge [black,thick]; [.{}  \edge [black,thick]; [.{}  \edge [black,thick]; [.{}  \edge [black,thick]; [.{}  \edge [black,thick]; [.{} ] \edge [gray!50]; [.{} ] \edge [gray!50]; [.{} ] \edge [gray!50]; [.{} ] \edge [gray!50]; [.{} ]] \edge [gray!50]; [.{}  \edge [gray!50]; [.{} ] \edge [gray!50]; [.{} ] \edge [gray!50]; [.{} ] \edge [gray!50]; [.{} ]] \edge [gray!50]; [.{}  \edge [gray!50]; [.{} ] \edge [gray!50]; [.{} ] \edge [gray!50]; [.{} ]] \edge [gray!50]; [.{}  \edge [gray!50]; [.{} ] \edge [gray!50]; [.{} ]] \edge [gray!50]; [.{} ]] \edge [black,thick]; [.{}  \edge [black,thick]; [.{}  \edge [black,thick]; [.{} ] \edge [black,thick]; [.{} ] \edge [gray!50]; [.{} ] \edge [gray!50]; [.{} ]] \edge [gray!50]; [.{}  \edge [gray!50]; [.{} ] \edge [gray!50]; [.{} ]] \edge [gray!50]; [.{}  \edge [gray!50]; [.{} ] \edge [gray!50]; [.{} ]] \edge [gray!50]; [.{} ]] \edge [gray!50]; [.{}  \edge [gray!50]; [.{}  \edge [gray!50]; [.{} ] \edge [gray!50]; [.{} ]] \edge [gray!50]; [.{}  \edge [gray!50]; [.{} ]] \edge [gray!50]; [.{} ]] \edge [gray!50]; [.{}  \edge [gray!50]; [.{} ]] \edge [gray!50]; [.{} ]] \edge [black,thick]; [.{}  \edge [black,thick]; [.{}  \edge [black,thick]; [.{}  \edge [black,thick]; [.{} ] \edge [gray!50]; [.{} ] \edge [gray!50]; [.{} ] \edge [gray!50]; [.{} ]] \edge [gray!50]; [.{}  \edge [gray!50]; [.{} ] \edge [gray!50]; [.{} ] \edge [gray!50]; [.{} ]] \edge [gray!50]; [.{}  \edge [gray!50]; [.{} ]] \edge [gray!50]; [.{} ]] \edge [gray!50]; [.{}  \edge [gray!50]; [.{}  \edge [gray!50]; [.{} ] \edge [gray!50]; [.{} ]] \edge [gray!50]; [.{}  \edge [gray!50]; [.{} ]] \edge [gray!50]; [.{} ]] \edge [gray!50]; [.{}  \edge [gray!50]; [.{} ]] \edge [gray!50]; [.{} ]] \edge [gray!50]; [.{}  \edge [gray!50]; [.{}  \edge [gray!50]; [.{}  \edge [gray!50]; [.{} ]] \edge [gray!50]; [.{} ]] \edge [gray!50]; [.{}  \edge [gray!50]; [.{} ]] \edge [gray!50]; [.{} ]] \edge [gray!50]; [.{}  \edge [gray!50]; [.{} ]] \edge [gray!50]; [.{} ]] \edge [black,thick]; [.{}  \edge [black,thick]; [.{}  \edge [black,thick]; [.{}  \edge [black,thick]; [.{}  \edge [black,thick]; [.{} ] \edge [gray!50]; [.{} ] \edge [gray!50]; [.{} ]] \edge [gray!50]; [.{}  \edge [gray!50]; [.{} ]] \edge [gray!50]; [.{} ]] \edge [gray!50]; [.{}  \edge [gray!50]; [.{}  \edge [gray!50]; [.{} ]]] \edge [gray!50]; [.{} ]] \edge [black,thick]; [.{}  \edge [black,thick]; [.{}  \edge [black,thick]; [.{}  \edge [black,thick]; [.{} ] \edge [gray!50]; [.{} ] \edge [gray!50]; [.{} ]] \edge [gray!50]; [.{}  \edge [gray!50]; [.{} ]] \edge [gray!50]; [.{} ]] \edge [gray!50]; [.{}  \edge [gray!50]; [.{} ]] \edge [gray!50]; [.{} ]] \edge [gray!50]; [.{}  \edge [gray!50]; [.{} ]] \edge [gray!50]; [.{} ]] \edge [gray!50]; [.{}  \edge [gray!50]; [.{}  \edge [gray!50]; [.{}  \edge [gray!50]; [.{} ]]] \edge [gray!50]; [.{} ]] \edge [gray!50]; [.{}  \edge [gray!50]; [.{} ]] \edge [gray!50]; [.{} ]] \edge [gray!50]; [.{}  \edge [gray!50]; [.{}  \edge [gray!50]; [.{}  \edge [gray!50]; [.{}  \edge [gray!50]; [.{}  \edge [gray!50]; [.{} ] \edge [gray!50]; [.{} ]] \edge [gray!50]; [.{}  \edge [gray!50]; [.{} ]] \edge [gray!50]; [.{} ]] \edge [gray!50]; [.{}  \edge [gray!50]; [.{} ]] \edge [gray!50]; [.{} ]] \edge [gray!50]; [.{}  \edge [gray!50]; [.{} ]] \edge [gray!50]; [.{} ]] \edge [gray!50]; [.{}  \edge [gray!50]; [.{}  \edge [gray!50]; [.{} ]]] \edge [gray!50]; [.{} ]] \edge [gray!50]; [.{}  \edge [gray!50]; [.{}  \edge [gray!50]; [.{}  \edge [gray!50]; [.{} ]]] \edge [gray!50]; [.{} ]] \edge [gray!50]; [.{}  \edge [gray!50]; [.{} ]] \edge [gray!50]; [.{} ]] \edge [black,thick]; [.{}  \edge [black,thick]; [.{}  \edge [black,thick]; [.{}  \edge [black,thick]; [.{}  \edge [black,thick]; [.{}  \edge [black,thick]; [.{}  \edge [black,thick]; [.{} ] \edge [gray!50]; [.{} ] \edge [gray!50]; [.{} ] \edge [gray!50]; [.{} ]] \edge [gray!50]; [.{}  \edge [gray!50]; [.{} ] \edge [gray!50]; [.{} ] \edge [gray!50]; [.{} ]] \edge [gray!50]; [.{}  \edge [gray!50]; [.{} ] \edge [gray!50]; [.{} ]] \edge [gray!50]; [.{}  \edge [gray!50]; [.{} ]]] \edge [gray!50]; [.{}  \edge [gray!50]; [.{}  \edge [gray!50]; [.{} ] \edge [gray!50]; [.{} ] \edge [gray!50]; [.{} ]] \edge [gray!50]; [.{}  \edge [gray!50]; [.{} ] \edge [gray!50]; [.{} ]] \edge [gray!50]; [.{} ]] \edge [gray!50]; [.{}  \edge [gray!50]; [.{}  \edge [gray!50]; [.{} ] \edge [gray!50]; [.{} ]] \edge [gray!50]; [.{} ]] \edge [gray!50]; [.{}  \edge [gray!50]; [.{} ]]] \edge [gray!50]; [.{}  \edge [gray!50]; [.{}  \edge [gray!50]; [.{}  \edge [gray!50]; [.{} ] \edge [gray!50]; [.{} ]] \edge [gray!50]; [.{}  \edge [gray!50]; [.{} ] \edge [gray!50]; [.{} ]] \edge [gray!50]; [.{} ]] \edge [gray!50]; [.{}  \edge [gray!50]; [.{} ]] \edge [gray!50]; [.{} ]] \edge [gray!50]; [.{}  \edge [gray!50]; [.{}  \edge [gray!50]; [.{} ]] \edge [gray!50]; [.{} ]] \edge [gray!50]; [.{} ]] \edge [gray!50]; [.{}  \edge [gray!50]; [.{}  \edge [gray!50]; [.{}  \edge [gray!50]; [.{}  \edge [gray!50]; [.{} ] \edge [gray!50]; [.{} ]] \edge [gray!50]; [.{} ]] \edge [gray!50]; [.{}  \edge [gray!50]; [.{} ]]] \edge [gray!50]; [.{}  \edge [gray!50]; [.{}  \edge [gray!50]; [.{} ]] \edge [gray!50]; [.{} ]] \edge [gray!50]; [.{} ]] \edge [gray!50]; [.{}  \edge [gray!50]; [.{} ]] \edge [gray!50]; [.{} ]] \edge [gray!50]; [.{}  \edge [gray!50]; [.{}  \edge [gray!50]; [.{}  \edge [gray!50]; [.{}  \edge [gray!50]; [.{}  \edge [gray!50]; [.{} ]] \edge [gray!50]; [.{} ]] \edge [gray!50]; [.{} ]] \edge [gray!50]; [.{} ]] \edge [gray!50]; [.{}  \edge [gray!50]; [.{} ]] \edge [gray!50]; [.{} ]] \edge [gray!50]; [.{}  \edge [gray!50]; [.{} ]] \edge [gray!50]; [.{} ]] \edge [gray!50]; [.{}  \edge [gray!50]; [.{}  \edge [gray!50]; [.{}  \edge [gray!50]; [.{}  \edge [gray!50]; [.{}  \edge [gray!50]; [.{} ]]] \edge [gray!50]; [.{} ]] \edge [gray!50]; [.{} ]] \edge [gray!50]; [.{} ]] \edge [gray!50]; [.{}  \edge [gray!50]; [.{} ]] \edge [gray!50]; [.{} ]] \edge [black,thick]; [.{}  \edge [black,thick]; [.{}  \edge [black,thick]; [.{}  \edge [black,thick]; [.{}  \edge [black,thick]; [.{}  \edge [black,thick]; [.{}  \edge [black,thick]; [.{}  \edge [black,thick]; [.{} ] \edge [black,thick]; [.{} ] \edge [gray!50]; [.{} ]] \edge [gray!50]; [.{}  \edge [gray!50]; [.{} ] \edge [gray!50]; [.{} ]] \edge [gray!50]; [.{}  \edge [gray!50]; [.{} ]]] \edge [gray!50]; [.{}  \edge [gray!50]; [.{}  \edge [gray!50]; [.{} ] \edge [gray!50]; [.{} ]] \edge [gray!50]; [.{}  \edge [gray!50]; [.{} ]]] \edge [gray!50]; [.{}  \edge [gray!50]; [.{} ]]] \edge [black,thick]; [.{}  \edge [black,thick]; [.{}  \edge [black,thick]; [.{}  \edge [black,thick]; [.{} ] \edge [gray!50]; [.{} ]] \edge [gray!50]; [.{}  \edge [gray!50]; [.{} ]]] \edge [gray!50]; [.{}  \edge [gray!50]; [.{} ]]] \edge [gray!50]; [.{}  \edge [gray!50]; [.{} ]]] \edge [gray!50]; [.{}  \edge [gray!50]; [.{}  \edge [gray!50]; [.{}  \edge [gray!50]; [.{}  \edge [gray!50]; [.{} ]]] \edge [gray!50]; [.{}  \edge [gray!50]; [.{} ]]] \edge [gray!50]; [.{} ]] \edge [gray!50]; [.{}  \edge [gray!50]; [.{} ]]] \edge [gray!50]; [.{}  \edge [gray!50]; [.{}  \edge [gray!50]; [.{}  \edge [gray!50]; [.{} ]]] \edge [gray!50]; [.{} ]] \edge [gray!50]; [.{} ]] \edge [black,thick]; [.{}  \edge [black,thick]; [.{}  \edge [black,thick]; [.{}  \edge [black,thick]; [.{}  \edge [black,thick]; [.{}  \edge [black,thick]; [.{}  \edge [black,thick]; [.{} ] \edge [gray!50]; [.{} ]] \edge [gray!50]; [.{} ]] \edge [gray!50]; [.{} ]] \edge [gray!50]; [.{} ]] \edge [gray!50]; [.{} ]] \edge [gray!50]; [.{} ]] \edge [gray!50]; [.{} ]] \edge [gray!50]; [.{}  \edge [gray!50]; [.{}  \edge [gray!50]; [.{} ]]] \edge [gray!50]; [.{} ]] \edge [black,thick]; [.{}  \edge [black,thick]; [.{}  \edge [black,thick]; [.{}  \edge [black,thick]; [.{}  \edge [black,thick]; [.{}  \edge [black,thick]; [.{}  \edge [black,thick]; [.{}  \edge [black,thick]; [.{}  \edge [black,thick]; [.{} ] \edge [gray!50]; [.{} ]] \edge [gray!50]; [.{}  \edge [gray!50]; [.{} ]]] \edge [gray!50]; [.{}  \edge [gray!50]; [.{}  \edge [gray!50]; [.{} ]]]] \edge [gray!50]; [.{}  \edge [gray!50]; [.{}  \edge [gray!50]; [.{} ]]]] \edge [gray!50]; [.{}  \edge [gray!50]; [.{}  \edge [gray!50]; [.{} ]]]] \edge [gray!50]; [.{}  \edge [gray!50]; [.{} ]]] \edge [gray!50]; [.{}  \edge [gray!50]; [.{} ]]] \edge [gray!50]; [.{} ]] \edge [gray!50]; [.{} ]] \edge [black,thick]; [.{}  \edge [black,thick]; [.{}  \edge [black,thick]; [.{}  \edge [black,thick]; [.{}  \edge [black,thick]; [.{}  \edge [black,thick]; [.{}  \edge [black,thick]; [.{}  \edge [black,thick]; [.{}  \edge [black,thick]; [.{}  \edge [black,thick]; [.{} ]]]]]]]]]]]]\end{tikzpicture}}

  \caption{Tree structure of the semigroup tree up to genus $11$.}
  \label{fig:tree}
  \end{figure}
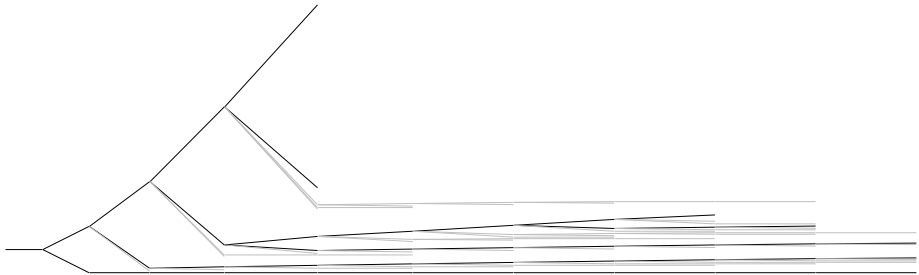

An infinite chain in the semigroup tree is an infinite sequence of subsequently included numerical semigroups that are adjacent in the semigroup tree. They have been analyzed in \cite{BB2009,RBsubmitted}.
Clearly, a numerical semigroup has infinitely many descendants in the
semigroup tree if and only if it lies in an infinite chain.
Computation shows that semigroups in infinite chains are rare, although the tree grows exponentially.
In Figure~\ref{fig:tree} one can see the tree structure of the semigroup tree up to genus 11. The edges corresponding to infinithe chains have been highlighted. It can be appreciated that they constitute indeed, a very small portion of the whole set of edges up to this genus.
It is the object of this work to prove that, for each fixed genus $g$, the portion of numerical semigroups of genus $g$ in infinite chains among the whole set of numerical semigroups of genus $g$ approaches $0$ as the genus grows to infinite. This means that most numerical semigroups have a finite number of descendants in the semigroup tree.

The multiplicity of a numerical semigroup is its smallest nonzero element.
Recent references analyze the expected multiplicity in terms of the genus of a semigroup as the genus grows to infinity, as well as the expected Frobenius number in terms of the multiplicity \cite{KaplanYe,Singhal}. Indeed, let $\gamma = (5 + \sqrt{5})/10$.
Fix $\varepsilon>0$. It is stated in \cite[Proposition 16]{KaplanYe}
that the portion of semigroups of genus $g$ with multiplicity between $(\gamma-\varepsilon)g$ and $(\gamma+\varepsilon)g$
with respect to the whole set of numerical semigroups of genus $g$ approaches $1$ as the genus $g$ grows to infinity. On the other hand, \cite[Theorem 4(1)]{KaplanYe} tells that the portion of semigroups of genus $g$ with Frobenius number $F$ satisfying $(2-\varepsilon)m<F<(2+\varepsilon)m$, where $m$ is the multiplicity of the semigroup, with respect to the whole set of numerical semigroups of genus $g$ approaches $1$ as the genus $g$ grows to infinity. See also \cite[Theorem 8]{Singhal}.

The left elements of a numerical semigroup are all the elements of the semigroup that are smaller than the Frobenius number.
It is proved in \cite{BB2009,RBsubmitted} that a numerical semigroup belongs to an ininite chain if and only if its nonzero left elements are not coprime.
In our proof we will use the aforementioned asymptotic results on the multiplicity and the Frobenius number of a numerical semigroup together with this characterization of semigroups in infinite chains.

\section{Proof of the rarity of numerical semigroups in infinite chains}

Let $\gamma = (5 + \sqrt{5})/10$.
Let $S_g$ be the set of numerical semigoups of genus $g$ and let $n_g=\#S_g$.
For a numerical semigroup $S$ let $m(S)$, $F(S)$, $L(S)$ denote its multiplicity, Frobenius number, and set of left elements, respectively.

Define

\begin{eqnarray*}
  A_g^m(\varepsilon_1)&=&\{S\in S_g:(\gamma-\varepsilon_1)g<m(S)<(\gamma+\varepsilon_1)g\}\\
  A_g^F(\varepsilon_2)&=&\{S\in S_g:(2-\varepsilon_2)m(S)<F(S)<(2+\varepsilon_2)m(S)\}\\
  A_g^{m,F}(\varepsilon_1,\varepsilon_2)&=&A_g^m(\varepsilon_1)\cap A_g^F(\varepsilon_2)\\
  S_g^{\bar\infty}&=&\{S\in S_g:\gcd(L(S))=1\}\\
  B_g(\varepsilon)&=&\{S\in S_g:(\gamma-\varepsilon)g<m(S)<(\gamma+\varepsilon)g\mbox{ and }(2-\varepsilon)\gamma g<F(S)<(2+\varepsilon)\gamma g\},
\end{eqnarray*}

\begin{lemma}{\cite[Proposition 16]{KaplanYe}}\label{l:KYm}
For any fixed $\varepsilon_1>0$, $\lim_{g\to\infty}\frac{\#A_g^m(\varepsilon_1)}{n_g}=1$.
\end{lemma}

\begin{lemma}{\cite[Theorem 4(1)]{KaplanYe}}\label{l:KYF}
For any fixed $\varepsilon_2>0$, $\lim_{g\to\infty}\frac{\#A_g^F(\varepsilon_2)}{n_g}=1$.
\end{lemma}

\begin{lemma}\label{l:limit}
For any fixed pair $\varepsilon_1>0,\varepsilon_2>0$, $\lim_{g\to\infty}\frac{\#A_g^{m,F}(\varepsilon_1,\varepsilon_2)}{n_g}=1$.
\end{lemma}

\begin{proof}
  It is a consequence of the fact that
  $$S_g=A^{m,F}_g(\varepsilon_1,\varepsilon_2)  \sqcup (A^m_g(\varepsilon_1)\setminus A^F_g(\varepsilon_2)) \sqcup (A^F_g(\varepsilon_2)\setminus A^m_g(\varepsilon_1))\sqcup(S_g\setminus((A^m_g(\varepsilon_1)\cup A^F_g(\varepsilon_2)))$$ and so,
  $$\#A^{m,F}_g(\varepsilon_1,\varepsilon_2)=n_g
  -\#(A^m_g(\varepsilon_1)\setminus A^F_g(\varepsilon_2))
  -\#(A^F_g(\varepsilon_2)\setminus A^m_g(\varepsilon_1))
  -\#(S_g\setminus(A^m_g(\varepsilon_1)\cup A^F_g(\varepsilon_2))),$$
  and the fact that, by Lemma~\ref{l:KYm} and Lemma~\ref{l:KYF},
  $$\lim_{g\to\infty}\frac{
    \#(A^m_g(\varepsilon_1)\setminus A^F_g(\varepsilon_2))}{n_g}\leq\lim_{g\to\infty}\frac{\#(S_g\setminus A^F_g(\varepsilon_2))}{n_g}=0$$

  $$\lim_{g\to\infty}\frac{\#(A^F_g(\varepsilon_2)\setminus A^m_g(\varepsilon_1))}{n_g}\leq\lim_{g\to\infty}\frac{\#(S_g\setminus A^m_g(\varepsilon_1))}{n_g}=0$$

  $$\lim_{g\to\infty}\frac{\#(S_g\setminus(A^m_g(\varepsilon_1)\cup A^F_g(\varepsilon_2)))}{n_g}\leq \lim_{g\to\infty}\frac{\#(S_g\setminus(A^m_g(\varepsilon_1))}{n_g}=0$$
\end{proof}

\begin{lemma}\label{l:dosepsilonsunepsilon}
  For any fixed $\varepsilon>0$, setting $\varepsilon_1=\frac{\gamma\varepsilon}{4+\varepsilon}$ and $\varepsilon_2=\frac{\varepsilon}{2}$, it holds
  $$A_g^{m,F}(\varepsilon_1,\varepsilon_2)  \subseteq B_g(\varepsilon).$$
\end{lemma}

\begin{proof}
  On one hand, $\varepsilon_1<\varepsilon$, and so, it is easy to check that any numerical semigroup $S$ in $A_g^{m,F}(\varepsilon_1,\varepsilon_2)$ satisfies $(\gamma-\varepsilon)<m(S)<(\gamma+\varepsilon)$.
  On the other hand,
  \begin{eqnarray*}F(S)&>&(2-\varepsilon_2)(\gamma-\varepsilon_1)g
    \\&=&\left(2-\frac{\varepsilon}{2}\right)\left(\gamma-\frac{\gamma\varepsilon}{4+\varepsilon}\right)g
    \\&=&2\gamma g-\frac{\gamma\varepsilon g}{2}-\frac{2\gamma\varepsilon g}{4+\varepsilon}+\frac{\gamma\varepsilon^2g}{2(4+\varepsilon)}
    \\&=&2\gamma g-\frac{\gamma\varepsilon g}{2}-\frac{4\gamma\varepsilon g-\gamma\varepsilon^2g}{2(4+\varepsilon)}
    \\&=&2\gamma g-\gamma\varepsilon g+\left(\frac{\gamma\varepsilon g}{2}-\frac{\gamma\varepsilon g(4-\varepsilon)}{2(4+\varepsilon)}\right)
    \\&>&(2-\varepsilon)\gamma g.\end{eqnarray*}
  while
  \begin{eqnarray*}F(S)&<&(2+\varepsilon_2)(\gamma+\varepsilon_1)g\\
    &=&\left(2+\frac{\varepsilon}{2}\right)\left(\gamma+\frac{\gamma\varepsilon}{4+\varepsilon}\right)g\\
    &=&2\gamma g+\frac{\gamma\varepsilon g}{2}+\frac{2\gamma\varepsilon g}{4+\varepsilon}+\frac{\gamma\varepsilon^2g}{2(4+\varepsilon)}
    \\&=&2\gamma g+\frac{\gamma\varepsilon g}{2}+\frac{4\gamma\varepsilon g+\gamma\varepsilon^2g}{2(4+\varepsilon)}
    \\&=&2\gamma g+\frac{\gamma\varepsilon g}{2}+\frac{\gamma\varepsilon g(4+\varepsilon)}{2(4+\varepsilon)}\\&=&(2+\varepsilon)\gamma g.\end{eqnarray*}
  
\end{proof}

\begin{lemma}\label{l:consecutiveelements}
  If $g\geq 6$ and $0<\varepsilon\leq \frac{3\gamma-2-\frac{1}{g}}{1+3\gamma}$, 
then any semigroup in $B_g(\varepsilon)$ has two consecutive elements.
\end{lemma}

\begin{proof}
First of all, notice that for $g\geq 6$ the expression $\frac{3\gamma-2-\frac{1}{g}}{1+3\gamma}$ is positive and, so, some value $\varepsilon$ exists in the interval of the hypothesis.
  
  Let $S\in B_g(\varepsilon)$. Let $L^*(S)=L(S)\setminus\{0\}$. It holds
  \begin{eqnarray*}\#L^*(S)&=&F(S)-g\\&>& ((2-\varepsilon)\gamma-1)g\\&=&(2\gamma-1-\varepsilon\gamma)g.\end{eqnarray*}
  If $S$ has no two consecutive elements, then \begin{eqnarray*}F(S)&\geq& m(S)+2\#L^*(S)-1\\&>&(\gamma-\varepsilon)g+2(2\gamma-1-\varepsilon\gamma)g-1\\
    &=&(5\gamma-\varepsilon-2-2\varepsilon\gamma)g-1\\
    &=&((2+\varepsilon)\gamma+3\gamma-2-\varepsilon(1+3\gamma))g-1\\
    &\geq&((2+\varepsilon)\gamma+3\gamma-2-(3\gamma-2-\frac{1}{g}))g-1\\
    &=&(2+\varepsilon)\gamma g.\end{eqnarray*}
  This contradicts the upper bound on the Frobenius number of the elements in $A_g(\varepsilon)$. Hence, $S$ must have two consecutive elements.
  \end{proof}

\begin{lemma}{\cite[Theorem 10]{BB2009},\cite[Lemma 3.4]{RBsubmitted}}\label{l:infinitechainsgcdL}
A numerical semigroup $S$ belongs to an infinite chain if and only if the elements in $L^*(S)$ are not coprime.
  \end{lemma}

\begin{lemma}\label{l:notinfinitechains}
  If $g\geq 6$ and $0<\varepsilon\leq \frac{3\gamma-2-\frac{1}{g}}{1+3\gamma}$, 
then $B_g(\varepsilon)\subseteq S_g^{\bar\infty}$.
\end{lemma}

\begin{proof}
It is a direct consequence of Lemma~\ref{l:consecutiveelements} and Lemma~\ref{l:infinitechainsgcdL}.
  \end{proof}

\begin{theorem}
$\lim_{g\to\infty}\frac{\#S_g^{\bar\infty}}{n_g}=1$.
\end{theorem}

\begin{proof}
  Fix $\varepsilon$ with $0<\varepsilon\leq \frac{3\gamma-2-\frac{1}{g}}{1+3\gamma}$. Set $\varepsilon_1=\frac{\gamma\varepsilon}{4+\varepsilon}$ and $\varepsilon_2=\frac{\varepsilon}{2}$.
  As a consequence of Lemma~\ref{l:dosepsilonsunepsilon} and Lemma~\ref{l:notinfinitechains},
  $$A_g^{m,F}(\varepsilon_1,\varepsilon_2)  \subseteq B_g(\varepsilon)\subseteq S_g^{\bar\infty},$$
  hence,
  $$\lim_{g\to\infty}\frac{\#S_g^{\bar\infty}}{n_g}\geq\lim_{g\to\infty}\frac{\#A_g^{m,F}(\varepsilon_1,\varepsilon_2)}{n_g}$$
  
  Now, by Lemma~\ref{l:limit},
  $$\lim_{g\to\infty}\frac{\#S_g^{\bar\infty}}{n_g}\geq\lim_{g\to\infty}\frac{\#A_g^{m,F}(\varepsilon_1,\varepsilon_2)}{n_g}=1,$$
  which concludes the proof.
  \end{proof}

\section*{Acknowledgment}
This work was supported by the Ministerio de Ciencia e Innovación, under grant PID2021-124928NB-I00, and by the AGAUR under grant 2021 SGR 00115. 
The second author was supported by the AGAUR under grant 2021 FISDUR 00189. 

The graph has been drawn using the drawsgtree tool, which can be downloaded from https://github.com/mbrasamoros/drawsgtree.

\bibliographystyle{plain}

\begin{thebibliography}{10}

\bibitem{BernardiniTorres}
Matheus Bernardini and Fernando Torres.
\newblock Counting numerical semigroups by genus and even gaps.
\newblock {\em Discrete Math.}, 340(12):2853--2863, 2017.

\bibitem{BlancoGarciaPuerto}
V\'{\i}ctor Blanco, Pedro~A. Garc\'{\i}a-S\'{a}nchez, and Justo Puerto.
\newblock Counting numerical semigroups with short generating functions.
\newblock {\em Internat. J. Algebra Comput.}, 21(7):1217--1235, 2011.

\bibitem{Br:Fibonacci}
M.~Bras-Amor{\'o}s.
\newblock Fibonacci-like behavior of the number of numerical semigroups of a
  given genus.
\newblock {\em Semigroup Forum}, 76(2):379--384, 2008.

\bibitem{Br:Bounds}
M.~Bras-Amor{\'o}s.
\newblock Bounds on the number of numerical semigroups of a given genus.
\newblock {\em J. Pure Appl. Algebra}, 213(6):997--1001, 2009.

\bibitem{seeds1}
M.~Bras-Amor{\'o}s and J.~Fern{\'a}ndez-Gonz{\'a}lez.
\newblock Computation of numerical semigroups by means of seeds.
\newblock {\em Math. Comp.}, 87(313):2539--2550, 2018.

\bibitem{rgd}
M.~Bras-Amor\'{o}s and J.~Fern\'{a}ndez-Gonz\'{a}lez.
\newblock The right-generators descendant of a numerical semigroup.
\newblock {\em Math. Comp.}, 89(324):2017--2030, 2020.

\bibitem{seeds2}
Maria Bras-Amor\'{o}s.
\newblock On the seeds and the great-grandchildren of a numerical semigroup.
\newblock {\em Math. Comp.}, 93(345):411--441, 2024.

\bibitem{BB2009}
Maria Bras-Amor{\'o}s and Stanislav Bulygin.
\newblock Towards a better understanding of the semigroup tree.
\newblock {\em Semigroup Forum}, 79(3):561--574, 2009.

\bibitem{DelgadoEliahouFromentin}
Manuel Delgado, Shalom Eliahou, and Jean Fromentin.
\newblock A verification of wilf's conjecture up to genus 100.
\newblock 2023.

\bibitem{EliahouFromentin:gapsets}
Shalom Eliahou and Jean Fromentin.
\newblock Gapsets and numerical semigroups.
\newblock {\em J. Combin. Theory Ser. A}, 169:105129, 19, 2020.

\bibitem{EliahouFromentin:gapsetsm}
Shalom Eliahou and Jean Fromentin.
\newblock Gapsets of small multiplicity.
\newblock In {\em Numerical semigroups}, volume~40 of {\em Springer INdAM
  Ser.}, pages 63--82. Springer, Cham, [2020] \copyright 2020.

\bibitem{EliahouRamirezAlfonsin}
Shalom Eliahou and Jorge Ram\'{\i}rez~Alfons\'{\i}n.
\newblock On the number of numerical semigroups {$\langle a,b\rangle$} of prime
  power genus.
\newblock {\em Semigroup Forum}, 87(1):171--186, 2013.

\bibitem{Elizalde}
S.~Elizalde.
\newblock Improved bounds on the number of numerical semigroups of a given
  genus.
\newblock {\em J. Pure Appl. Algebra}, 214(10):1862--1873, 2010.

\bibitem{FromentinHivert}
J.~Fromentin and F.~Hivert.
\newblock Exploring the tree of numerical semigroups.
\newblock {\em Math. Comp.}, 85(301):2553--2568, 2016.

\bibitem{Kaplan}
N.~Kaplan.
\newblock Counting numerical semigroups by genus and some cases of a question
  of {W}ilf.
\newblock {\em J. Pure Appl. Algebra}, 216(5):1016--1032, 2012.

\bibitem{KaplanYe}
Nathan Kaplan and Lynnelle Ye.
\newblock The proportion of {W}eierstrass semigroups.
\newblock {\em J. Algebra}, 373:377--391, 2013.

\bibitem{ODorney}
E.~O'Dorney.
\newblock Degree asymptotics of the numerical semigroup tree.
\newblock {\em Semigroup Forum}, 87(3):601--616, 2013.

\bibitem{RBsubmitted}
Mariana Rosas-Ribeiro and Maria Bras-Amor\'os.
\newblock Infinite chains in the tree of numerical semigroups.
\newblock Submitted.

\bibitem{Singhal}
Deepesh Singhal.
\newblock Distribution of genus among numerical semigroups with fixed
  {F}robenius number.
\newblock {\em Semigroup Forum}, 104(3):704--723, 2022.

\bibitem{Zhai}
A.~Zhai.
\newblock Fibonacci-like growth of numerical semigroups of a given genus.
\newblock {\em Semigroup Forum}, 86(3):634--662, 2013.

\bibitem{Zhao}
Y.~Zhao.
\newblock Constructing numerical semigroups of a given genus.
\newblock {\em Semigroup Forum}, 80(2):242--254, 2010.

\bibitem{Zhu}
Daniel~G. Zhu.
\newblock Sub-fibonacci behavior in numerical semigroup enumeration.
\newblock 2023.

\end{thebibliography}

\end{document}